\documentclass[9pt]{article}

\usepackage{amssymb, amsmath, latexsym, amscd}
\usepackage{pdfsync}
\usepackage{amsthm}  

\def\sref#1{Section~\ref{#1}}

\def\tref#1{Theorem~\ref{#1}}

\def\cref#1{Corollary~\ref{#1}}
\def\pref#1{Proposition~\ref{#1}}
\def\lref#1{Lemma~\ref{#1}}
\def\eref#1{Example~\ref{#1}}

\def\sref#1{Section~\ref{#1}}

\def\tref#1{Theorem~\ref{#1}}

\def\cref#1{Corollary~\ref{#1}}
\def\pref#1{Proposition~\ref{#1}}
\def\lref#1{Lemma~\ref{#1}}
\def\eref#1{Example~\ref{#1}}

\newtheoremstyle{thm}{1em}{1em}{\itshape}{}{\bfseries}{.}{.5em}{} 
\newtheoremstyle{rem}{1em}{1em}{}{}{\bfseries}{.}{.5em}{}           

\theoremstyle{thm}                
\newtheorem{thm}{Theorem}[section]
\newtheorem{prop}[thm]{Proposition}

\newtheorem{lem}[thm]{Lemma}

\newtheorem{dfn}[thm]{Definition}

\theoremstyle{rem} 
\newtheorem{rmk}[thm]{Remark}

\newtheorem{ex}[thm]{Example}

\newcommand{\intro}[1]
{\renewcommand{\thesection}{\fnsymbol{section}}
\setcounter{section}{-1}
\section{#1}
\renewcommand{\thesection}{\arabic{section}}
}

\DeclareMathOperator{\Ker}{Ker}
\DeclareMathOperator{\im}{Im}

\DeclareMathOperator{\Dim}{dim}

\DeclareMathOperator{\id}{id}

\newcommand{\Pf}{\noindent{\bf Proof}. }
\newcommand{\cqfd}
{%
\mbox{}%
\nolinebreak%
\hfill%
\rule{2mm}{2mm}%
\medbreak%
\par%
}
\renewcommand{\a}{{\mathcal A}{}} 
\renewcommand{\b}{{\mathcal B}{}}

\newcommand{\F}{{\mathcal F}{}}
\newcommand{\G}{{\mathcal G}{}}

\newcommand{\N}{{\mathbb N}{}}

\newcommand{\R}{\mathbb R}

\newcommand{\torus}{\mathbb T}
\newcommand{\U}{{\mathcal U}{}}

\newcommand{\ZE}{{\mathcal Z}{}}
\newcommand{\Z}{\mathbb Z}

\newcommand{\eps}{\varepsilon }
\newcommand{\ol}{\overline}
\newcommand{\ul}{\underline}
\newcommand{\hot}{\hat{\otimes}}

\newcommand{\Wlog}{without loss of generality }
\newcommand{\wrt}{with respect to }
\newcommand{\rp}{respectively }

\begin{document} 

\title{Some remarks on a K\"unneth formula for foliated de Rham cohomology}

\author{M\'elanie Bertelson$^\ast$} 

\date{April 8, 2008}

\maketitle

\begin{center}
$^\ast$ Chercheur Qualifi\'e F.N.R.S. \\
D\'epartement de Math\'ematiques, C.P. 218 \\
Universit\'e Libre de Bruxelles \\
 Boulevard du Triomphe \\
1050 Bruxelles \\
Belgique\\
{\tt mbertels@ulb.ac.be}
\end{center}

\begin{abstract}
The K\"unneth formula is one of the basic tools for computing cohomology. Its validity for foliated cohomology, that is, for the tangential de Rham cohomology of a foliated manifold, is investigated. The main difficulty encountered is the non-Hausdorff nature of the foliated cohomology spaces, forbidding  the completion of the tensor product. The results presented here are a K\"unneth formula when both factors have Hausdorff foliated cohomology, a K\"unneth formula when one factor has Hausdorff finite-dimensional foliated cohomology and a counterexample to an alternative version of the K\"unneth formula. The proof of the second result involves a right inverse for the foliated de Rham differential.
\end{abstract}

\noindent
\intro{Introduction}

The tangential de Rham cohomology or {\em foliated cohomology} of a foliated manifold $(M,\F)$ is the cohomology of the complex obtained by forming the quotient of the Fr\'echet space of ordinary smooth forms on the manifold by those who vanish along the leaves of the foliation. Our initial interest for this cohomology comes from the observation that its vanishing in degree two may, under certain circumstances, be an obstruction to existence of a foliated symplectic structure, or equivalently a regular Poisson structure whose underlying foliation is $\F$\footnote{cf.~Bertelson, M., Foliations associated to regular Poisson structures. {\em Commun. Contemp. Math.} {\bf 3}, No. 3 (2001) 441--456.}. Among the tools for computing de Rham cohomology is the K\"unneth formula which asserts that the cohomology space of a product is isomorphic to the completed tensor product of the cohomology spaces of the factors via the map
$$\varphi : \bigoplus_{p+q = n} H^p(M) \otimes H^q(N) \to H^n(M \times N) : a \otimes b \mapsto a \wedge b.$$
This map indeed induces a map on foliated cohomology but because these spaces do not generally enjoy the Hausdorff separation property, the completion of the tensor product may not even be defined. \\

The results obtained so far and exposed in the present paper are~:

\begin{enumerate}
\item[-] The K\"unneth formula is valid when the foliated cohomology spaces of both factors are Hausdorff. This is a consequence of a result due to Grothendieck and exposed in \cite{S}. We have nevertheless included a relatively detailed proof in \sref{K-H}.

\item[-] It is also valid when the foliated cohomology of one of the factors is finite-dimensional and Hausdorff. Notice that it is not necessary to complete the tensor product in that case. This result was already known when one of the factors is a one-leaf foliation (\cite{EK} or \cite{M&S}). Our proof requires the construction of a right inverse for the foliated de Rham differential. It is based on results in the theory of splitting of exact sequences of Fr\'echet spaces (\cite{M&V} and \cite{V}).

\item[-] In the simple case where one of the factors, say $(M,\F)$, has a non-Hausdorff foliated cohomology and the other factor, say $(N,\G)$, is a manifold foliated by its points, a natural alternative version to the K\"unneth formula would involve $C^\infty (N, H^*({\cal F}))$ in place of the completed tensor product. Nevertheless, we have constructed on the torus ${\mathbb T}^2$ foliated by Liouville slope lines a smooth family of exact forms --- that represents thus the zero element in $C^\infty (N, H^*({\cal F}))$ --- which is not the coboundary of any continuous family of functions --- that corresponds therefore to a non-zero element in $H^*({\cal F} \times \G)$. 
\end{enumerate}

Many relevant questions remain unanswered, among which~:

\begin{enumerate}
\item[-] Does a more sophisticated version of the K\"unneth formula, involving some type of higher order functors, hold in a non-Hausdorff situation~? 
\item[-] Does the foliated de Rham differential still admit a right inverse when the assumption of finite-dimensionality of the foliated cohomology is relaxed~? 
\item[-] Can a foliated manifold have a finite-dimensional non-Hausdorff cohomology or does finite-dimensional imply Hausdorff for foliated cohomology~?
\end{enumerate}

Finally, the results of this paper may apply to other cohomologies. We have in mind the Poisson cohomology of a Poisson manifold (not surprisingly, the tangential Poisson cohomology of a regular Poisson structure is isomorphic to the foliated cohomology of the induced foliation). For instance, the 
K\"unneth formula for Poisson cohomology is valid when the cohomology spaces are Hausdorff.
 
\section*{Acknowledgments} 

I wish to thank warmly Professor Dietmar Vogt for helping me to clarify how I could use his results on splitting of exact sequences of Fr\'echet spaces exposed in \cite{V} and \cite{M&V}. I also wish to thank Alan Weinstein and Pierre Bieliavsky for useful discussions regarding the subject.

\section{Preliminaries}

Let $(M, \F)$ be a foliated manifold, that is, a smooth Hausdorff second countable manifold $M$ endowed with a smooth foliation $\F$. The space of smooth $p$-forms, $p\geq 0$, is denoted by $\Omega^p(M)$ (a smooth $0$-form is just a smooth function) and the space of all forms by $\Omega^*(M)$. The weak $C^\infty$ topology provides $\Omega^p(M)$ (and $\Omega^*(M)$) with the structure of a nuclear Fr\'echet space. We are interested in the nuclear property because it guarantees uniqueness of the completion of the tensor product with any other Fr\'echet space.\\

Let us recall that a Fr\'echet space is a locally convex, metrizable, complete topological vector space. We will not attempt to explain the nuclear property here, but rather refer to Sections 47 and 50 in \cite{T}. For our purpose it is sufficient to know that the set of smooth functions on an open subset of ${\mathbb R}^n$ is nuclear (cf.~\cite{T}[Corollary of Theorem 51.5, p~530]), that a product of nuclear spaces is nuclear and that a Hausdorff projective limit of nuclear spaces is nuclear (cf.~\cite{T}[Proposition 50.1, p~514]). Indeed, $\Omega^p(M)$ is the projective limit of the spaces $\Omega^p(\phi_\alpha(U_\alpha))$, where $(U_\alpha, \phi_\alpha)$ runs through an atlas on $M$. We will occasionally  write TVS for topological vector space.\\

Consider the space $\Omega^p(M, \F) = \{\omega \in \Omega^p(M) ; \omega|_F \; \forall \mbox{ leaf } F \}$ of forms vanishing along the leaves of $\F$. It is a closed subspace of $\Omega^p(M)$. Thus the  quotient $\Omega^p(M) / \Omega^p(M, \F)$ is a Fr\'echet nuclear space as well (cf.~\cite{T}[p~85 and Proposition 50.1, p~514]). It is the space of {\em foliated $p$-forms}. The de Rham differential $d : \Omega^*(M) \to \Omega^{*+1}(M)$ which is a continuous linear map induces the {\em foliated de Rham differential} $d_\F : \Omega^*(\F) \to \Omega^{*+1}(\F)$ with like properties. The space of $d_\F$-closed (\rp $d_\F$-exact ) foliated $p$-forms is denoted by $\ZE^p(\F)$ (\rp $\b^p(\F)$). The cohomology $H^*(\F) =  \ZE^*(\F) / \b^*(\F)$ is called the {\em foliated (de Rham) cohomology} of $(M, \F)$.

\begin{rmk} The (ordinary) de Rham differential is always a homomorphism, that is, the image of an open subset of $\Omega^p(M)$ under $d$ consists of a relative open subset of $d(\Omega^p(M))$. This is a consequence of the fact that a form is exact if and only if its integral over any closed cycle vanishes, showing that exact forms are a closed subset which by the open mapping theorem for metrizable and complete topological vector spaces implies that $d$ is open (cf.~\cite{T}[Theorem 17.1, p~170]). In contrast, the differential $d_\F$ needs not be a homomorphism, as illustrated by the \eref{torus} of the Liouville slopes foliations on the torus $\torus^2$. Observe that assuming that $d_\F$ is open is equivalent to assuming that $\b^*(\F)$ is closed (by the open mapping theorem for one direction and the observation that the image by a homomorphism of a complete metrizable TVS is a closed 
space for the other direction) or that the cohomology $H^*(\F)$ is Hausdorff. 
\end{rmk}

With regards to the previous remark, the following examples are useful to keep in mind.
\begin{ex}(Kronecker foliations)\label{torus}
Consider the foliation of $\R^2$ by parallel lines of slope $\alpha \in \R$. Being invariant under the action of $\Z^2$ by translations, this foliation induces a foliation, denoted $\F_\alpha$, on the torus $\torus^2$. The leaves are circles when $\alpha$ is a rational number and are dense lines otherwise. 
The foliated de Rham cohomology of $\F_\alpha$ for $\alpha$ irrational depends on the type of irrational number considered. More specifically, it is infinite-dimensional and non-Hausdorff  (with a one-dimensional Hausdorff quotient) when $\alpha$ is a Liouville number or one-dimensional and Hausdorff otherwise. The proof  of that well-known fact can be found in \cite{H} or \cite{M&S} and will appear implicitly in \sref{ex}. Let us remind the reader of the definition of a Liouville number.
\end{ex}

\begin{dfn} A Liouville number $\alpha$ is an irrational number that is {\it well-ap\-pro\-xi\-ma\-ble} by rational numbers. More precisely, for all integers $p \geq 1$, there exist relatively prime integers $m, n,$ with $n > 1$ such that

$$\big| \alpha - \frac{m}{n} \big| < \frac{1}{|n|^p}$$

\end{dfn}
A typical example of such a number is Liouville's constant $\sum_{k=1}^\infty 10^{-k!}$. Liouville numbers are transcendental because an algebraic number $\alpha$ of degree $p \geq 2$ admits a constant $c$ such that 
$$\big| \alpha - \frac{m}{n} \big| > \frac{c}{|n|^p},$$
for all integers $m, n$ with $n>0$. On the other hand $e$ and $\pi$ for instance are not Liouville, as are uncountably many transcendental numbers. The set of Liouville numbers is a countable intersection of open dense sets and has measure zero. A non-Liouville number is sometimes called a {\it generic number}. 

\begin{ex}\label{vanishing-cycle} Let $(M, \F)$ be a foliation that has a vanishing $k$-cycle, that is, a  smooth foliated map $v : (S^k \times [0,1], \F_\pi) \to (M, \F)$, where $S^k$ is a sphere of dimension $k$ and $\F_\pi$ is the foliation by the fibers of the canonical projection $\pi : S^k \times [0,1] \to [0,1]$, such that the image of $S^k \times \{t\}$ is homotopically trivial in its leaf for each $t$ except $t=0$. A $p$-dimensional foliation from which a point is removed carries a vanishing $(p-1)$-cycle. We explain hereafter, in the specific case of a punctured foliation $(M, \F) = (N-\{q\}, \G|_{N-\{q\}})$, how the presence of the vanishing $(p-1)$-cycle implies that $H^p(\F)$ is non-Hausdorff and infinite-dimensional. The argument can certainly be extended to a larger class of vanishing cycles. \\

Observe that our vanishing cycle can be ``filled", in the sense that there exists a foliated map $\ol{v} : (D^p \times [0,1] - {\rm int}D^p \times \{0\}, \F_\pi) \to (M, \F)$ that extends $v$. Let $\Omega$ be a foliated volume form on $(N, \G)$ and let $f$ be a smooth function on $M$ approaching infinity near the puncture. Then $f \Omega$ is a foliated closed $p$-form on $M$ than cannot be foliated exact. Indeed, suppose on the contrary that $f \Omega = d_\F \alpha$. Then, by Stokes' theorem, 
$$\int_{\ol{v}(D^p \times \{t\})} f \Omega = \int_{v(S^{p-1} \times \{t\})} \alpha.$$
Clearly, as $t$ approaches $0$, the right-hand side of the previous equality converges to $\int_{v(S^{p-1} \times \{0\})} \alpha$ while the left-hand side diverges, yielding a contradiction. Besides, it is not too difficult to construct an example of a non-exact $p$-form of this type that is the limit of a sequence of exact forms, showing that the set of foliated exact forms is not closed in the set of foliated closed forms.
\end{ex}

\section{K\"unneth formula when the cohomology is Hausdorff}\label{K-H}

The main result of the present section, that is, a K\"unneth formula for foliated cohomology when the foliated cohomology of each factor is Hausdorff, is not original as it is essentially a consequence of a theorem due to Grothendieck and exposed in \cite{S}. (A proof in terms of sheaf can also be found in the literature, namely in Glen E.~Bredon, {\it Sheaf theory.} Second edition. Graduate Texts in Mathematics~170. Springer-Verlag, New York, 1997.) We have nevertheless decided to write here a relatively detailed explanation of it, with systematic references to the book \cite{T} for the background functional analysis, believing that some readers might find it useful to have the proof expressed in a language familiar to differential geometers with references from just one very well-written book. \\

Let $(M, \F)$ and $(N, \G)$ be two foliated manifolds both having the property that the foliated de Rham differential is a homomorphism. Consider the (algebraic) tensor product $\Omega^p(\F) \otimes \Omega^q(\G)$. There are two natural ways to construct a topology on the tensor product of two locally convex Hausdorff topological vector spaces, namely the $\eps$ and the $\pi$ topology (cf.~\cite{T}[Section 42 and 43]), thus yielding two different completions of the tensor product. However, when one of  the factors is Fr\'echet nuclear, both topologies coincide (cf.~\cite{T}[Theorem 50.1, p~511]). So in our case we can ignore this issue and write $\Omega^p(\F) \hot \Omega^q(\G)$ for the completion --- \wrt this unique natural topology --- of the  tensor product of $\Omega^p(\F)$ with $\Omega^q(\G)$. Moreover, the tensor product of two continuous linear maps $f_1:E_1 \to F_1$ and $f_2 : E_2 \to F_2$ between nuclear Fr\'echet spaces is a continuous linear map $f_1 \otimes f_2 : E_1 \otimes E_2 \to F_1 \otimes F_2 \subset F_1 \hot F_2$ (\cite{T}[Proposition 43.6, p~439]) which induces a continuous linear map $f_1\hot f_2 : E_1 \hot E_2 \to F_1 \hot F_2$ between the completions. In general, the completion of a Hausdorff locally convex TVS $E$ is denoted by $\hat{E}$ and the extension of a continuous linear map $u : E \to F$ to the completions by $\hat{u} : \hat{E} \to \hat{F}$ 
(\cite{T}[Theorem 5.1, p~39]).\\

Consider the tensor product complex $(\Omega^*(\F)\hot \Omega^*(\G), d)$ defined as follows~:

$$\bigl(\Omega^*(\F) \hot \Omega^*(\G)\bigr)^n \stackrel{\rm def}{=} \bigoplus_{p+q = n} 
\Omega^p(\F) \hot \Omega^q(\G),$$
with differential $d = d_\F \hot 1 + \eps \hot d_\G$, where $\eps (\omega) = (-1)^p \omega$ when 
$\omega$ is a foliated form of degree $p$. It follows from general considerations that  $\Omega^*(\F) \hot \Omega^*(\G)$ is a nuclear Fr\'echet space (cf.~\cite{T}[Proposition 50.1 p~514) as well and that $d$ is a homomorphism. The latter assertion is a consequence of \cite{T}[Proposition 43.9, p~441] and the fact that the sum of two homomorphisms is a homomorphism. \\

There is a natural map $\varphi$ between and $\Omega^*(\F) \hot \Omega^*(\G)$ and $\Omega^*(\F \times \G)$, given by extension of the map
$$\begin{array}{ccccc}
\ul{\varphi} & : & \Omega^*(\F) \otimes \Omega^*(\G) & \longrightarrow & \Omega^*(\F \times \G) 
\\
&& \displaystyle{\sum_{i=1}^I } \alpha_i \otimes \beta_i & \longrightarrow & {p_M}^*(\alpha_i) \wedge {p_N}^*(\beta_i),
\end{array}$$
where $p_M$ and $p_N$ denote the projections of $M \times N$ onto $M$ and $N$ respectively. It is clearly a cochain map, that is, $\varphi \circ d = d_{\F \times \G} \circ \varphi$ and therefore induces a map on foliated cohomology.

\begin{thm}\label{Kunneth}(K\"unneth formula) The map $\varphi$ is an isomorphism on cohomology~: 
$$H^n(\F \times \G) \cong \bigr( H^*(\F) \hot H^*(\G) \bigl)^n.$$
\end{thm}
This is a direct consequence of the following two results~:

\begin{thm}\label{Grothendieck}(Grothendieck) \cite{S} Let $(E^*, d_E)$, $(F^*, d_F)$ be two differential complexes of Fr\'echet spaces and homomorphisms. Suppose that the $E^p$'s are nuclear. Consider the differential complex $(E^* \hot F^*, d)$ with $d = d_E \hot 1 + \eps \hot d_F$. Then $H^*(E \hot F) \cong H^*(E) \hot H^*(F)$.
\end{thm}

\begin{prop}\label{isocompl} The differential complexes $\bigr(\Omega^*(\F)\hot \Omega^*(\G), d \bigl)$ and $\bigr( \Omega^*(\F \times \G), d_{\F \times \G} \bigl)$ are isomorphic under the map $\varphi$.
\end{prop}

The proof of \tref{Grothendieck} relies mostly on the next two lemmas~:

\begin{lem}\label{tpses} Let $E$, $F$, $G$ and $H$ be four Fr\'echet spaces with either $E$, $F$ and $G$ nuclear or $H$ nuclear. Let $u : E \to F$ and $v : F \to G$ be linear homomorphisms such that 
$0 \to E \stackrel{u}{\to} F \stackrel{v}{\to} G \to 0$ is a short exact sequence. Then the sequence
$$0 \to E \hot H \stackrel{u \hot \id}{\longrightarrow} F\hot H \stackrel{v \hot \id}{\longrightarrow} G \hot H \to 0$$
is a short exact sequence of Fr\'echet spaces and linear homomorphisms as well. 
\end{lem}

\Pf The fact that $u \hot \id$ (\rp $v \hot \id$) is a 1-1 (\rp onto) map follows from \cite{T}[Proposition 43.6, p~440] (\rp \cite{T}[Proposition 43.9, p~441]). Exactness at $F \hot H$ is argued as follows. Firstly, observe that 
$$0 \to E \otimes H \stackrel{u \otimes \id}{\longrightarrow} F \otimes H \stackrel{v \otimes \id}{\longrightarrow} G \otimes H \to 0$$
is a short exact sequence of homomorphisms. Indeed, the corollary of Proposition 43.7, p~441 in \cite{T} implies that $u \otimes \id$ is a homomorphism. As for $v \otimes \id$, it suffices to know that a basis of neighborhoods of $0$ for the  $\pi$-topology consists of the convex hulls of sets of type $U \otimes V = \{u \otimes v ; u \in U \mbox{ and } v \in V\}$, where $U$ (\rp $V$) is a balanced neighborhood of $0$ in the first factor (\rp second factor), that is, $U$ is a neighborhood of $0$ such that $\lambda u \in U, \; \forall \; |\lambda| \leq 1, u \in U$. Therefore, $G \otimes H \cong F \otimes H / u \otimes \id (E \otimes H)$. Secondly, it is not difficult to prove that if $E$ is a metrizable TVS and if $N \subset E$ is a closed subspace then 
$$\widehat{E/N} \cong \hat{E} / \hat{N},$$ where $\hat{N}$ denotes the closure of $N$ in the completion $\hat{E}$ of $E$.
\cqfd

\begin{lem}\label{les-ass-to-ses} Let $(A^*, d_A)$, $(B^*, d_B)$ and $(C^*, d_C)$ be three differential complexes of metrizable complete TVS's and homomorphisms and let

$$0 \to A^* \stackrel{f}{\to} B^* \stackrel{g}{\to} C^* \to 0$$
be  a short exact sequence of differential complexes with $f$, $g$ continuous maps (hence homomorphisms by the open mapping theorem). Then the usual long exact sequence 

$$... \to H^*(A) \stackrel{f_*}{\to} H^*(B) \stackrel{g_*}{\to} H^*(C) \stackrel{\nu}{\to} H^{*+1}(A) \to ...$$
is well-defined with $f_*$, $g_*$ and $\nu$ homomorphisms.
\end{lem}

\Pf Since $d_A$, $d_B$ and $d_C$ are homomorphisms, all spaces involved (i.e.~cocycles, coboundary and quotients of the formers by the latters) are complete metrizable spaces. The open mapping theorem implies thus that any surjective continuous linear map between those spaces will be a homomorphism. The only things that requires a proof is therefore the continuity of $\nu$ which is easily seen by chasing open sets in the diagram providing the construction of $\nu$, namely, 

$$\begin{array}{ccccccccc}
0 & \to & A^{p+1} & \stackrel{f}{\to} & B^{p+1} & \stackrel{g}{\to} & C^{p+1} & \to & 0 \\
&& \uparrow d_A && \uparrow d_B && \uparrow d_C && \\
0 & \to & A^p & \stackrel{f}{\to} & B^p & \stackrel{g}{\to} & C^p & \to & 0. \\
\end{array}$$
\cqfd

\noindent
{\bf Proof of \tref{Grothendieck}} Let us introduce some notation~:
$$\begin{array}{ccccccc}
\ZE^p = E^p \cap \Ker d_E &&& \b^p = E^p \cap \im d_E &&& H^p = \ZE^p /\b^p \\
\ZE'^p = F^p \cap \Ker d_F&&& \b'^p = F^p \cap \im d_F &&& H'^p = \ZE'^p /\b'^p 
\end{array}$$
Now consider the following exact sequences of linear homomorphisms~:
$$\begin{array}{lcr}
0 \to \ZE^p \to E^p \to \b^{p+1} \to 0  & \qquad &  0 \to \ZE'^p \to F^p \to \b'^{p+1} \to 0 \\
0 \to \b^p \to \ZE^p \to H^p \to 0 & \qquad & 0 \to \b'^p \to \ZE'^p \to H'^p \to 0
\end{array}$$
By \lref{tpses}, they induce the following other exact sequences of linear homomorphisms (obtained by tensoring with a fixed space and the identity map), where we have omitted the superscripts $^*$~:
$$\begin{array}{llllllllllr}
0 & \to & (\ZE \hot F)^n & \to & (E \hot F)^n & \to & (\b \hot F)^{n+1} & \to & 0 && (1)\label{3}\\
0 & \to & (\ZE \hot \ZE')^n & \to & (\ZE \hot F)^n & \to & (\ZE \hot \b')^{n+1} & \to & 0 && (2) \\
0 & \to & (\ZE \hot \b')^n & \to & (\ZE \hot \ZE')^n & \to & (\ZE \hot H')^n & \to & 0 && (3)\\
0 & \to & (\b \hot \ZE')^n & \to & (\b \hot F)^n & \to & (\b \hot \b')^{n+1} & \to & 0 && (4)\\
0 & \to & (\b \hot \b')^n & \to & (\b \hot \ZE')^n & \to & (\b \hot H')^n & \to & 0 && (5)\\
0 & \to & (\b \hot H')^n & \to & (\ZE \hot H')^n & \to & (H \hot H')^n & \to & 0 && (6) 
\end{array}$$

The first one is also an exact sequence of differential complexes when $(\ZE \hot F)^*$ (\rp $(\b \hot 
F)^*$) is endowed with $d' = \eps \hot d_F$ (\rp $d'' = - \eps \hot d_F$), yielding, by \lref{les-ass-to-ses}, the long exact sequence 
\setcounter{equation}{6}
\begin{equation}\label{les}
... \to H^*(\ZE \hot F) \to H^*(E \hot F) \to H^{*+1}(\b \hot F) \to H^{*+1}(\ZE \hot F) \to ...
\end{equation}
Moreover, the sequences (2) and (3) imply that $H^*(\ZE \hot F) \cong (\ZE \hot H')^*$. Indeed the sequence (2) identifies $(\ZE \hot \ZE')^*$ (\rp $(\ZE \hot \b')^*$) as being the kernel (\rp the image) of the differential $d'$ (the $\eps$ does not affect that conclusion since all maps are graded). Moreover, sequence (3) says that the quotient of $(\ZE \hot \ZE')^*$ by $(\ZE \hot \b')^*$ is isomorphic to $(\ZE \hot H')^*$. Similarly (4) and (5) imply that $H^*(\b \hot F) \cong (\b \hot H')^*$. With these isomorphisms, the sequence (\ref{les}) becomes
\begin{equation}\label{lesbis}
... \to (\ZE \hot H')^* \to H^*(E \hot F) \to (\b \hot H')^{*+1} \stackrel{\nu}{\to} (\ZE \hot H')^{*+1} \to ...
\end{equation}
We will prove that $\nu$ is the map induced by the natural inclusion $\b^* \to \ZE^*$. Indeed, consider the following diagram~:

$$\begin{array}{ccccccccc}
0 & \to & (\ZE \hot F)^{n+1} & \to & (E \hot F)^{n+1} & \to & (\b \hot F)^{n+1} & \to & 0 \\
&& \uparrow d' && \uparrow d && \uparrow d'' && \\
0 & \to & (\ZE \hot F)^n & \to & (E \hot F)^n & \to & (\b \hot F)^n & \to & 0. \\
\end{array}$$
Pick $\sum_{i = 1}^k b_i \otimes z_i$ in $\b^p \otimes \ZE'^q$. Let $b_i = d_E x_i$ for some $x_i$ in $E^{p-1}$, then $d(\sum_{i = 1}^k x_i \otimes z_i) = \sum_{i = 1}^k b_i \otimes z_i$. This shows that $\nu$ and $i \hot \id$ coincide on the subspace $\b^p \otimes H'^q$. Therefore they coincide on all of $\b \hot H'$. \\

Since $\nu$ is an injective map (again by \cite{T}[Proposition 43.7, p~440]), the long exact sequence (\ref{lesbis}) is equivalent to the short exact sequence

$$0 \to (\b \hot H')^* \to (\ZE \hot H')^* \to H^*(E \hot F) \to 0.$$
Hence $H^*(E \hot F) \cong (\ZE \hot H')^* / (\b \hot H')^*$ and the latter space is isomorphic to $H \hot H'$, as the sequence (6) shows.
\cqfd

\noindent
{\em Proof of \pref{isocompl}} The proof in notationally heavy but conceptually quite simple. First observe that the continuous map $\ul{\varphi} : \Omega^*(\F) \otimes \Omega^*(\G) \to \Omega^*(\F \times \G)$ is injective. We will prove hereafter that it is a homomorphism with dense image, implying that its extension $\varphi : \Omega^*(\F) \hot \Omega^*(\G) \to \Omega^*(\F \times \G)$ is an isomorphism. \\

To prove that $\varphi$ is a homomorphism, recall that the following subsets of $\Omega^p(\F)$ form a basis of neighborhoods of $0$~:
$$\begin{array}{rl}
\U (r, \eps, \{(U_i, \phi_i)\}, \{K_i\}) = & \{\omega \in \Omega^p(\F); |D^a \omega_{i, j_1 ... j_p} (x)| \leq \eps \\ 
& \forall \mbox{ multi-index } a = (a_1, ..., a_{\Dim M})\mbox{ with } |a| \leq r, \\
& \forall \; 1 \leq i \leq n, \forall \; x \in K_i\},
\end{array}$$
where $r$ is some non-negative integer, where $\eps > 0$, where $\{(U_i, \phi_i); 1 \leq i \leq n\}$ is a  finite collection of foliated charts and where $K_i$ is a compact subset of $U_i$ for each $1 \leq i \leq n$. The functions $\omega_{i, j_1 ... j_p}$, $1 \leq j_1< ...< j_p \leq \Dim M$ denote the tangential coordinates of $\omega$ \wrt the chart $(U_i, \phi_i)$ and $D^a \omega_{i, j_1 ... j_p}$ is the $a^{\rm th}$ derivative
$$D^a \omega_{i, j_1 ... j_p} = \frac{\partial^{a_1}}{\partial x_1} ... \frac{\partial^{a_{\dim M}}}{\partial x_{\dim M}} (\omega_{i, j_1 ... j_p}).$$ 
We want to verify that if $U$ is a neighborhood of $0$ in $\Omega^p(\F) \otimes \Omega^q(\G)$, then $\varphi(U) \supset O \cap \varphi (\Omega^p(\F) \otimes \Omega^q(\G))$, for some neighborhood of $0$ in $\Omega^{p+q}(\F \times \G)$. Now a neighborhood of $0$ in $\Omega^p(\F) \otimes \Omega^q(\G)$ can be chosen of the type 
$$\U(U^o, V^o) = \Bigl\{\sum_{i=1}^I \alpha_i \otimes \beta_i ; \sup_{x' \in U^o, y' \in V^o} \Bigl|\sum_{i=1}^I <x', \alpha_i><y', \beta_i> \Bigr| \leq 1 \Bigr\},$$
where $U$ (\rp $V$) is a neighborhood of $0$ in $\Omega^p(\F)$ (\rp $\Omega^q(\G)$) and where $U^o$ denotes the polar of $U$, that is, the following subset of the dual $\Omega^p(\F)'$ of $\Omega^p(\F)$~: $U^o = \{x' \in \Omega^p(\F)' ; |<x', u>| \leq 1 \; \forall \; u \in U \}$ (same for $V^o$). Let us say a few words about this issue in the next paragraph. \\

The (algebraic) tensor product $E \otimes F$ of two locally convex Hausdorff TVS's is isomorphic to $B(E'_\sigma, F'_\sigma)$, the vector space of continuous bilinear forms on the product $E'_\sigma \times F'_\sigma$ of the duals of $E$ and $F$ each endowed with its respective weak topology (topology of pointwise convergence) (cf.~\cite{T}[Proposition 42.4, p~432]). The latter space can be naturally realized as a subspace of a complete locally convex Hausdorff TVS, namely the space $\b_\eps(E'_\sigma, F'_\sigma)$ of {\em separately continuous} bilinear forms on $E'_\sigma \times F'_\sigma$ with the $\eps$-topology, or topology of uniform convergence on products of equicontinuous subsets of $E'$ and $F'$ (cf.~\cite{T} Definition 43.1, p~434, Proposition 42.3, p~430). When endowed with the topology induced by $\b_\eps(E'_\sigma, F'_\sigma)$, the tensor product of $E$ and $F$ is denoted by 
$$E \otimes_\eps F.$$
 The topology on $\b_\eps(E'_\sigma, F'_\sigma)$ can be defined by the following basis of neighborhoods of $0$~:
$$\U (A, B) = \{\phi \in \b_\eps(E'_\sigma, F'_\sigma) ; |\phi(A, B)| \leq 1\},$$
where $A$ (\rp $B$) is an equicontinuous subset of $E'$ (\rp $F'$). The reason for the restriction to equicontinuous sets (rather than just bounded sets) is explained in \cite{T}[p~427--428]. Now any equicontinuous subset of $E'$ is contained in the polar $U^o$ of some neighborhood $U$ of $0$ (cf.~\cite{T}[Proposition 32.7, p~341]). Thus, a basis of neighborhoods of $0$ is also given by the sets $\U_\delta (U^o, V^o)$, where $U$ (\rp $V$) runs through a basis of neighborhoods of $0$ in $E$ (\rp $F$).

Returning to our proof that $\varphi$ is a homomorphism, we claim that if $U = \U (r, \eps, \{(U_i, \phi_i)\}, \{K_i\})$ and $V = \U (s, \delta, \{(V_k, \psi_k)\}, \{L_k\}) $, then $\varphi(\U(U^o, V^o)) \supset \im \varphi \cap O$, where $O = \U(\max\{r, s\}, \eps \delta, \{(U_i \times V_k, \phi_i \times \psi_k)\}, \{K_i \times L_k\})$. For that purpose, it will be useful to observe that the set $\U(r, \eps, \{(U_i, \phi_i)\}, \{K_i\})$ is the polar of the following subset of the dual of $\Omega^p(\F)$~: 
\begin{multline*}
\a (r, \eps, \{(U_i, \phi_i)\}, \{K_i\}) = \{\ell_{a, \eps, i, j_1 ...  j_p, x} (\omega) = \frac{1}{\eps} \partial^a \omega_{i, j_1 ...j_p} (x) ; \\ |a| \leq r, i = 1, ..., n, 1 \leq j_1< ... < j_p \leq \dim M, x \in K_i\}.
\end{multline*}
Thus $U = A^o$. On the other hand, a locally convex Hausdorff TVS $E$ is isomorphic to the dual of its weak dual, that is, $E \cong (E'_\sigma)'$ (cf.~\cite{T}[Proposition 35.1, p~361]), and if $U = A^o$, then $U^o = (A^o)^o$ coincides with the closed convex balanced hull of $A$ (i.e.~the closure of the convex hull of $\cup_{\{\lambda ; |\lambda| \leq 1\}} \lambda A$), denoted by $\Gamma A$ (cf.~\cite{T}[Proposition 35.3, p~362]). Furthermore, one verifies directly from the definitions involved that $\U(A,B) = \U(\Gamma A, \Gamma B)$. Thus, $\U(U^o, V^o) = \U(A, B)$ with $A = \a(r, \eps, \{(U_i, \phi_i)\}, \{K_i\})$ and $B = \a(s, \delta, \{(V_k, \psi_k)\}, \{L_k\})$. \\

Now let $\theta = \sum_{t=1}^T \alpha_t \otimes \beta_t \in \Omega^p(\F) \otimes \Omega^q(\G)$ be such that $\varphi(\theta)$ belongs to $O$, that is, $\forall$ multi-index $c$ with $|c| \leq rs$, $\forall \; i,  k, j_1< ... < j_p, l_1 < ... < l_q$, $\forall \; (x, y) \in K_i \times L_k$, one has

$$\Bigl| D^c \bigl( \sum_{t=1}^T (\alpha_t)_{i, j_1 ... j_p} (\beta_t)_{k, l_1 ... l_q}\bigr) (x,y) \Bigr| \leq \eps \delta, \;\; \mbox{or}$$
$$\Bigl| \sum_{t=1}^T \frac{1}{\eps} D^a (\alpha_t)_{i, j_1 ... j_p} (x) \frac{1}{\delta} D^b (\beta_t)_{k, l_1 ... l_q} (y) \Bigr| \leq 1,$$
where $c = (a_1, ..., a_{\Dim M}, b_1, ..., b_{\Dim N})$. Equivalently, 
$$\Bigl| \sum_{t=1}^T \bigl< \ell_{a, \eps, i, j_1 ... j_p, x}, \alpha_t \bigr> \bigl< \ell_{b, \delta, k, l_1 ... l_q, y}, \beta_t \bigl> \Bigr| \leq 1,$$
which means that $\theta \in \U(A, B)$, thus proving that $\varphi$ is a homomorphism. \\

It remains to prove that the image of $\varphi$ is dense in $\Omega^n(\F \times \G)$. It is essentially a consequence of the fact that polynomial functions are dense in the space of smooth functions on the euclidean space, implying that if $X$ and $Y$ are open subsets of $\R^n$ and $\R^m$ \rp then the tensor product of the spaces of smooth functions with compact supports on $X$ and $Y$, $C^\infty_c(X) \otimes C^\infty_c(Y)$, is dense in the space $C^\infty(X \times Y)$ of smooth functions on $X \times Y$ (cf.~\cite{T}[Theorem 39.2, p~409 and Corollary~1, p~159]). 
Let $\omega \in \Omega^{n}(\F \times \G)$ and consider $U$ a neighborhood of $\omega$ of the type $\omega + \U(r, \eps, \{(U_i \times V_k, \phi_i \times \psi_k), \{K_i \times L_k\}\})$. For each tangential component $\omega_{i,k,j_1 ... j_p, l_1 ... l_q}$, denoted hereafter $\omega_{i,k,J,L}$, of $\omega$ \wrt the chart $\phi_i \times \psi_k$, pick functions  $f_{i,k,J,L}^n \in C^\infty_c(U_i)$ and $g_{i,k,J,L}^n \in C^\infty_c(V_k)$, $n = 1, ..., N$ such that 
$$\sum_{n=1}^N f_{i,k,J,L}^n \; g_{i,k,J,L}^n$$ lies in $\omega_{i,k,J,L} + \U(r, \eps, K_i \times 
L_k) \subset C^\infty(U_i \times V_k)$. 
Then the form 
$$\sum_{n,J,L} f_{i,k,J,L}^n \; g_{i,k,J,L}^n dx^J \wedge dx^L$$
belongs to $U \cap \varphi(\Omega^p (\F) \otimes \Omega^q(\G))$. 
\cqfd

\section{K\"unneth formula when one of the factors is finite-dimensional and Hausdorff}

Another natural question is whether the K\"unneth formula holds in the case where the tensor product does not need to be completed, that is, when one of the factors, say $H^*(\F)$, is finite-dimensional. The answer is positive provided that factor is also Hausdorff. There is no assumption on the second factor. This statement was already well-known when $\F$ is a one-leaf foliation (cf.~\cite{EK} or \cite{M&S}). We use the fact that under the previous assumptions, the foliated de Rham differential $d_{\F}$ admits a right inverse, which is implied by results in the theory of splitting of exact sequences of Fr\'echet spaces appearing in \cite{M&V} and \cite{V}. As to the hypothesis that $H^*(\F)$ is Hausdorff, it suggests another  question~: ``Can $H^*(\F)$
 be finite-dimensional and non-Hausdorff?".

\begin{prop}\label{finiteK} Let $(M,\F)$ and $(N,\G)$ be foliated manifolds. Suppose that $H^*(\F)$ is finite-dimensional and Hausdorff. Then, as TVS's, 
$$H^*(\F \times \G) \cong H^*(\F) \otimes H^*(\G).$$
\end{prop}

\noindent
One may verify directly that, although we may not anymore quote theorems about coincidence of the $\eps$ and $\pi$-topology on $H^*(\F) \otimes H^*(\G)$ since $H^*(\F)$ is not Hausdorff, both these topologies  coincide with the direct sum topology that appears when $H^*(\F) \otimes H^*(\G)$ is identified with a finite direct sum $\oplus_{i=1}^n H^*(\G)$ via the choice of a basis of $H^*(\F)$.
 \\

\Pf The idea is to replace the complex $(\Omega^*(\F), d_{\F})$ by a homotopy equivalent finite-dimensional complex $(V, d_V)$. It is then easy to prove that $H^*(V \; \otimes \; \Omega^*(\G))$ coincides with $H^*(\Omega^*(\F) \; \hot \; \Omega^*(\G)) = H^*(\F \times \G)$. Besides, it is well-known that $H^*(V \otimes \Omega^*(\G)) = H^*(V) \otimes H^*(\Omega^*(\G))$ as vector spaces and it is not difficult to be convinced that this equality holds for the topologies as well. So we are done. To obtain an equivalence with a finite-dimensional complex we need a right inverse for the foliated differential $d_{\F}$, that is, a continuous linear map $\varphi : \b^{*+1}(\F) \to \Omega^*(\F)$ such that $d_{\F} \circ \varphi = \id$. This is the content of the \lref{right-inv} below. Let us assume this fact and proceed with the present proof. \\

The complex $(\Omega^*(\F), d_\F)$ is denoted hereafter by $(\Omega, d)$, $\Ker d_\F$ by $\ZE$ and $\im d_\F$ by $\b$. Consider closed foliated forms (of pure degree) $\alpha_1, ..., \alpha_n$ representing a basis  $\{[\alpha_1], ..., [\alpha_n]\}$ of $H^*(\F)$. The subset $V = \{\alpha_1, ..., \alpha_n\}$ endowed with the zero differential ($d_V = 0$) is a finite-dimensional subcomplex of $(\Omega, d)$ with cohomology $H^*(\F)$. It is thus (algebraically) homotopy equivalent to $(\Omega, d)$ (cf.~\cite{Sp}[Theorem 7.4.10, p~192]). We show hereafter that the homotopy, its inverse and the equivalence may be chosen continuous when $d$ admits a right inverse.\\

We first need to setup some notation.

\begin{enumerate}
\item[-] The natural inclusion $V \to \Omega$ is denoted by $i$.
\item[-] $\varphi : \b \to \Omega$ denotes a continuous linear right inverse to $d$.
\item[-] The cohomology class of a closed form $\beta$ is denoted by $[\beta]$.
\item[-] Since $H^*(\F)$ is Hausdorff and finite-dimensional, the linear map $e : H^*(\F) \to V$ such that $e([\alpha_i]) = \alpha_i$ is continuous; it is a right inverse for the natural projection ${\mathcal Z}^*(\F) \to H^*(\F)$ with values in $V$.
\end{enumerate}

Define $\sigma : \Omega \to V$ and $D: \Omega \to \Omega$ by 

\begin{enumerate}
\item[] $\sigma (\beta) = e [\beta -  \varphi(d \beta)]$.
\item[] $D(\beta) = -\varphi \Bigl( (\beta - \varphi(d \beta)) - i \circ e[ \beta - \varphi(d \beta)] \Bigr)$.
\end{enumerate}
The maps $\sigma$ and $D$ are clearly continuous. It is only necessary to verify that $\sigma$ is a cochain map, that $\sigma \circ i = \id_V$ and that $i \circ \sigma = \id_{\Omega} + D \circ d + d \circ D$. The first two assertions are obvious and the third one is proved hereafter.

$$\begin{array}{l}
\Bigl(D \circ d + d \circ D\Bigr) (\beta) \\
\begin{array}{lllllll}
&&&&& = & - \varphi (d \beta ) - d \circ \varphi 
\Bigl( (\beta - \varphi(d \beta)) - i \circ e [\beta - \varphi(d \beta)] \Bigr) \\
&&&&& = & -\beta + i \circ e [\beta - \varphi(d \beta)] \\
&&&&& = & - \beta + i \circ \sigma (\beta)
\end{array}
\end{array}$$

Now the continuous cochain maps $i$ and $\sigma$ induce continuous cochain maps $i \hot \id : V \hot \Omega^*(\G) \to \Omega^*(\F) \hot \Omega^*(\G)$ and $\sigma \hot \id : \Omega^*(\F) \hot \Omega^*(\G) \to V \hot \Omega^*(\G)$ such that 
$$(\sigma \hot \id) \circ (i \hot \id) = \id$$
and 
$$(i \hot \id) \circ (\sigma \hot \id) = \id + D' \circ d_{\F \times \G} + d_{\F \times \G} \circ D',$$ 
where $D' = D \hot \id$ and $d_{\F \times \G} = d_\F \hot \id + \eps \hot d_\G$. (The last assertion follows from the fact that $\eps \circ D + D \circ \eps = 0$.) Furthermore, the continuous cochain maps $i \hot \id$ and $\sigma \hot \id$ induce continuous maps on cohomology that are inverse to one another. This shows that $H^*(V \hot \Omega^*(\G)) \cong H^*(\F \times \G)$ as TVS's. Of course, $V \hot 
\Omega^*(\G)$ is the same as $V \otimes \Omega^*(\G)$ since $V \otimes E$ is complete when $V$ is finite-dimensional and $E$ is complete. \\

Finally, the fact that $H^*(V \otimes \Omega^*(\G)) \cong V \otimes H^*(\G)$ follows from considering the short exact sequences

$$0 \longrightarrow V \otimes \ZE^*(\G) \longrightarrow V \otimes \Omega^*(\G) \stackrel{\eps \otimes d_\G}{\longrightarrow} V \otimes \b^{*+1}(\G) \longrightarrow 0$$
and 
$$0 \longrightarrow V \otimes \b^*(\G) \longrightarrow V \otimes \ZE^*(\G) \longrightarrow V \otimes H^*(\G) \longrightarrow 0.$$
Observe that $\eps \otimes d_\G$ is continuous but not open. Likewise $V \otimes \b^{*+1}(\G)$ is not complete and $V \otimes H^*(\G)$ is not Hausdorff. Nevertheless, the first sequence tells us that the kernel of the differential $\eps \otimes d_\G$ on $V \otimes \Omega^*(\G)$ is $V \otimes \ZE^*(\G)$ and that its image is $V \otimes \b^{*+1}(\G)$. Besides, the topology induced on $V \otimes \ZE^*(\G)$ (\rp $V \otimes \b^{*+1}(\G)$) from its embedding in $V \otimes \Omega^*(\G)$ coincides with the tensor product topology (cf.~\cite{T}[Proposition 43.7, p 440]). Finally, since the maps in the second sequence are homomorphisms (remembering that $V \otimes H^*(\G)$ carries the direct sum topology), the quotient $V \otimes \ZE^*(\G) / V \otimes \b^*(\G)$ is isomorphic to $V \otimes H^*(\G)$.
\cqfd


Regarding existence of a right inverse for $d_\F$, we treat the simpler case of a compact underlying manifold first.

\begin{lem}\label{right-inv} If M is compact and $H^*(\F)$ is finite-dimensional and Hausdorff then $d_{\F}$ admits a continuous linear right inverse.
\end{lem}

\Pf The proof relies on the following result.

\begin{thm}[\cite{M&V} Splitting theorem 30.1, p~378]\label{splittingcpt} Let $E$, $F$, $G$ be Fr\'echet-Hilbert spaces and let $0 \to F \to G \to E \to 0$ be a short exact sequence of continuous linear maps. If $E$ has the property
(DN) and $F$ has the property $(\Omega)$, then the sequence splits.
\end{thm} 
\noindent
We explain hereafter why \tref{splittingcpt} can be applied to the short exact sequence 
$$0 \to {\mathcal Z}^*(\F) \to \Omega^*(\F) \to {\mathcal B}^*(\F) \to 0.$$ 

The assumption that the spaces are Fr\'echet-Hilbert is automatically satisfied for nuclear Fr\'echet spaces (see \cite{M&V}[Definition p~370 and Lemma 28.1 p~344]). We mention the definitions of properties (DN) and ($\Omega$) for completeness but we will only need here the fact that they are stable under performing certain operations. Let $E$ be a Fr\'echet space endowed with a countable fundamental systems $\{|| \cdot ||_k ; k \geq 1\}$ of seminorms (i.e.~for all $x \in E$, $x \neq 0$, there exists a $k$ such that $||x||_k > 0$ and for all $k_1, k_2$, there exist a $k_3$ and a constant $C$ such that $\max \{|| \cdot ||_{k_1}, || \cdot ||_{k_2}\} \leq C || \cdot ||_{k_3}$). The property (DN) is satisfied by $E$ if and only if it supports a continuous norm $|| \cdot ||$ on $E$ such that for any seminorm $|| \cdot ||_k$ there exists a constant $C$ and a seminorm $|| \cdot ||_K$ such that 
$$|| x ||_k^2 \leq C || x ||\; || x ||_K  \quad \forall x \in E.$$
The property ($\Omega$) is satisfied by $E$ if and only if for each $p \geq 1$ there exists a $q\geq 1$ so that for every $k \geq 1$, there exists a $0 < \theta < 1$ and a constant $C$ such that 
$$||y||^*_q \leq C {||y||^*_p}^{1 - \theta} {||y||^*_k}^{\theta} \quad \forall y \in E',$$
where $||y||^*_k$ means $\sup\{|y(x)| ; ||x||_k \leq 1\}$.\\

Both properties are satisfied by the Schwartz space 
$$s = \{(x_j)_{j\geq 1} ; \sum_{j=1}^\infty |x_j|^2 j^{2k} < \infty \; \forall k \geq 1\}$$ 
(\cite{M&V}[Example 29.5(1), p~363, Lemma 29.2, p~359 and Lemma 29.11, p~368). Besides, any space of $C^\infty$-sections of a finite-dimensional vector bundle $E$ over a compact manifold is isomorphic, as topological vector space, to $s$ (cf.~\cite{Va})\footnote{
The reference \cite{Va} contains a proof of the fact that for a compact manifold $M$, the space $C^\infty(M) \cong s$ which can easily be adjusted to the case of $C^\infty(M, E)$. Indeed, a finite partition of unity $\{\theta_i ; i = 1, ..., n\}$ subordinated to a cover of $M$ by trivializing open subsets allows us to identify the space of smooth sections of the bundle $E$ with a finite direct sum 
$\oplus_i C_c^\infty(C_i, {\mathbb R}^d)$, where $C_i$ is the support of $\theta_i$,  where $d$ is the rank of $E$ and where $C_c^\infty(C_i, {\mathbb R}^d)$ is the set of smooth functions with compact support in $C_i$. Because each $C_c^\infty(C_i, \mathbb R)$ is isomorphic to $s$ (\cite{Va}[(5) p 536]) and $s \oplus s \cong s$ (\cite{Va}[(5) p 327]), we reach our conclusion. 
}.

Thus the space of foliated forms $\Omega^*(\F)$ enjoys the properties (DN) and $(\Omega)$. Besides, property (DN) is inherited by closed subspaces (cf.~\cite{M&V}[Lemma 29.2, p~359]). So ${\mathcal B}^*(\F)$ has property (DN). To see that ${\mathcal Z}^{*}(\F)$ has property $(\Omega)$, we use the fact that $H^{*}(\F)$ is finite-dimensional. Indeed, the property $(\Omega)$ is inherited by quotients by closed subspaces (cf.~\cite{M&V}[Lemma 29.11(2), p~368]) so that ${\mathcal B}^{*}(\F)$, which is isomorphic to $\Omega^{*-1}(\F)/{\mathcal Z}^{*-1}(\F)$, has property $(\Omega)$. Since $H^{*}(\F)$ is finite-dimensional, the natural projection ${\mathcal Z}^{*}(\F) \to H^{*}(\F)$ admits a right inverse so that ${\mathcal Z}^{*}(\F) \cong {\mathcal B}^{*}(\F) \oplus H^{*}(\F)$ and can be thought of as a quotient of $\Omega^{*-1}(\F) \oplus H^{*}(\F)$ (by ${\mathcal Z}^{*-1}(\F) \oplus \{0\}$), which is itself also isomorphic to the Schwartz space~$s$ when $* \geq 1$. Finally, ${\mathcal Z}^0(\F)$ has property $(\Omega)$ because it is finite-dimensional (it is thus a Banach space).
\cqfd

When the manifold is not compact, the space $\Omega^*(\F)$ is isomorphic to $s^\N$ (argument similar to the compact case with a locally finite partition of unity subordinated to an open cover of $M$ by foliated chart domains. See also (11) p 438 in \cite{Va} for a proof that $C^\infty(M) \cong s^\N$) rather than $s$. One can nevertheless obtain a similar conclusion provided a condition is imposed on the homomorphism $A$, as the following statement shows.

\begin{thm}[\cite{V}, Theorem 3.5, p~820] Let $0 \to F \to G \to E \to 0$ be an exact sequence of nuclear Fr\'echet spaces, $A$ an $SK$-homomorphism. If $E$ has property (DN$_{\rm loc}$) and $F$ property $(\Omega_{loc})$, then the sequence splits.
\end{thm} 

A SK-homomorphism between two locally convex topological vector spaces $E$ and $F$ is a continuous linear map $A : E \to F$ such that for any seminorm $p$ on $E$, there is a seminorm $q$ on $F$ such that 
$$A (\Ker p) \supset \Ker q \cap \im A.$$ 

The condition on $A$ is in fact the weakest possible. Indeed, a continuous linear map between locally convex topological vector spaces that admits a continuous right inverse is necessarily a  SK-homomorphism (cf.~\cite{V}[Lemma 1.5, p~814]). 

\begin{rmk} The following example illustrates the fact that a surjective homomorphism is not necessarily a SK-homomorphism (cf.~\cite{V}[p~814]). Consider the continuous map 
$$C^\infty([-1,1]) \to \omega = \{(x_i)_{i \geq 0}\} : f \mapsto (\frac{d^if}{dt}(0))_{i \geq 0}.$$
It is clearly continuous and open for the usual Fr\'echet topologies but it is not a SK-homomorphism since $C^\infty([-1,1])$ admits continuous norms while $\omega$ does not. 
\end{rmk}

It might still be true that for the specific case of foliated forms, a homomorphism is automatically a SK-homomorphism. And it is so when the foliated de Rham cohomology is finite-dimensional, as asserted by the following result.
\begin{prop}
Let $(M,\F)$ be a foliated manifold whose foliated cohomology $H^*(\F)$ is finite-dimensional. Then, if the foliated de Rham differential $d_\F$ is a homomorphism, it is also a SK-homomorphism.
\end{prop}

\Pf For short, the space $\Omega^*(\F)$ (\rp $\b^*(\F)$) is denoted hereafter by $E$ (\rp $F$) and the differential $d_\F$ by $d$. We begin with the following preliminary observation. If $p$ is a continuous seminorm on $E$, its pushforward  $\ul{p}$, defined by 
$$\ul{p}(y) = \inf \{p(x) ; d(x) = y\},$$
is a continuous seminorm on $\im d = F$ (it is continuous if and only if $d$ is open). Its kernel is the set
$$\{y ; \mbox{ there exists a sequence } (x_k)_{k \geq 1} \subset d^{-1}(y) \mbox{ with } \lim_{k \to \infty} 
p(x_k) = 0 \}.$$ 
The inclusion $d(\Ker p) \subset \Ker \ul{p}$ is always true, unlike its reverse. \\

It is also useful to notice that it is sufficient to verify the condition characterizing a SK-homomorphism on the elements of a fundamental system of continuous seminorms on $E$\footnote{
If $p$ is any seminorm, there exists a $p_\alpha$ such that $U_p \supset \eps U_{p_\alpha}$, where $U_p = \{x \in E ; p(x) < 1\}$ (same for $U_{p_\alpha}$). In other words, $p \leq \frac{1}{\eps} p_\alpha$, which implies that $\Ker p \supset \Ker p_\alpha$. So if $A(\Ker p_\alpha) \supset \Ker q \cap \im A$ for some seminorm $q$ on $F$, then $A(\Ker p) \supset \Ker q \cap \im A$.}. We use the fundamental system of seminorms  described hereafter. Let 
$$p (\alpha) = p_{r, \{(U_i, \phi_i)\}, \{K_i\}} (\alpha) = \sup \{D^a \alpha_{i, j_1 ... j_p} (x); |a| \leq r, x\in 
K_i\},$$
where $r$ is some non-negative integer and $a$ is a multi-index $a = (a_1, ..., a_{\Dim M})$, where 
$\{(U_i, \phi_i); 1 \leq i \leq n\}$ is a  finite collection of foliated charts and $K_i$ is a compact subset of $U_i$ for each $1 \leq i \leq n$ and where the functions $\omega_{i, j_1 ... j_p}$, $1 \leq j_1< ...< j_p \leq \Dim M$ denote the tangential coordinates of $\omega$ \wrt the chart $(U_i, \phi_i)$. We assume that the compact set $\cup_i K_i$, denoted hereafter by $M_p$, is a submanifold with boundary of maximal dimension. To obtain a fundamental system of such seminorms, consider an exhaustion $M^1 \subset M^2 \subset ... $ of $M$ by compact submanifolds-with-boundary. Then cover each $M^\ell$ with finitely many foliated chart domains $U^\ell_1, ..., U^\ell_n$ and decompose $M^\ell$ into a finite union $M^\ell = K^\ell_1 \cup ... \cup K^\ell_n$ of compact sets such that $K^\ell_i \subset  U^\ell_i$. This procedure ensures that $\Ker p = \{\alpha ; \alpha|_{M^\ell} = 0\}$ when $p = p_{r, \{(U^\ell_i, \phi^\ell_i)\}, \{K^\ell_i\}}$ and that the set of foliated forms on $M_p$ can be defined easily. \\

Now let $p = p_{r, \{(U_i, \phi_i)\}, \{K_i\}}$. We claim that 
$$\Ker\ul{p} \subset d(\Ker p).$$
To prove this, let $\beta$ be an element of $\Ker\ul{p}$. This means that there exists a sequence $(\alpha_k)$ in 
$E$ with $d \alpha_k = \beta$ and $p(\alpha_k) \longrightarrow 0$ as $k \to \infty$. Intuitively, we are saying that the linear spaces $d^{-1}(\beta)$ and $\Ker p$ are in some sense asymptotic to one another. Our finite-dimensional intuition of linearity tells us that $d^{-1}(\beta)$ and $\Ker p$ should intersect. This would indeed be true if $\Ker d$ was finite-dimensional, which is far from being the case in general. Nevertheless, after forming the quotient by the set of exact forms we reach a finite-dimensional situation to which this observation can be applied.\\

Let $B = \b^*(\F) \subset \Ker d \subset E$ be the set of foliated exact forms. The quotient map $\pi : E \to E/B : x \to \ol{x}$ is a surjective linear homomorphism, as is the induced map $\ol{d} : E/B \to F : \ol{x} \to d \, x$. Define the seminorm $\ol{p} (\ol{x}) = \inf \{p(x') ; \pi(x') = \ol{x}\}$. Then the sequence $(\ol{\alpha}_k)$ satisfies
\begin{enumerate}
\item[-] $\ol{d}\ol{\alpha}_k = \beta$ for all $k$'s,
\item[-] $\lim_{k \to \infty}\ol{p} (\ol{\alpha}_k) = 0$.
\end{enumerate}
Let $\ol{\alpha}_o \in \ol{d}^{-1}(\beta)$. Then $\ol{d}^{-1}(\beta) = \ol{\alpha}_o + \Ker \ol{d}$ and thus $\ol{\alpha}_k = \ol{\alpha}_o + \tau_k$ for some $\tau_k \in \Ker \ol{d} = \pi (\Ker d) = H^*(\F)$. We consider two cases. \\

\noindent 
{\bf First case~:} $\Ker \ol{p} \cap \Ker \ol{d} = \{0\}$. \\

Since $\Ker \ol{d}$ is finite-dimensional, there exists a compact $K \subset \Ker \ol{d}$ not containing $0$ and such that any half line $\{\lambda \ol{\alpha} ; \lambda \in [0, \infty) \}$ intersects $K$ non-trivially (a sphere for some norm on $\Ker\ol{d}$ for instance). Let $\tau_k = a_k \rho_k$, where $a_k \geq 0$ and $\rho_k \in K$. Then $\rho_k$ admits a converging subsequence, say $\rho_k \to \rho$. If 
$\{a_k\}$ is bounded, then we are done. Otherwise, a subsequence of $(a_k)$, that we denote $(a_k)$ as well, converges to $\infty$ which yields a contradiction~:
$$\infty = \infty \, \ol{p} (\rho) \longleftarrow a_k \, \ol{p}(\rho_k) = \ol{p} (a_k \rho_k) \leq \ol{p} (\ol{\alpha}_o + a_k \rho_k) + \ol{p}(\ol{\alpha}_o) \longrightarrow \ol{p}(\ol{\alpha}_o).$$
So, the sequence $\ol{\alpha}_k$ admits a converging subsequence. Its limit is thus an element $\ol{\alpha}$ such that 
\begin{enumerate}
\item[-] $\ol{d}\ol{\alpha} = \beta$ and 
\item[-] $\ol{p}(\ol{\alpha}) = 0$.
\end{enumerate}

\noindent 
{\bf Second case~:} $\Ker \ol{p} \cap \Ker \ol{d} \neq \{0\}$. \\

Set $\Ker \ol{p} \cap \Ker \ol{d} = D_o$ and choose a complementary subspace $\Ker \ol{d} = D_o \oplus D_1$. Let $\tau_k = \eta_k + \mu_k$ with $\eta_k \in D_o$ and $\mu_k \in D_1$. Then $\ol{p}(\ol{\alpha}_o + \mu_k) \to 0$ as well. Indeed, 
$$\ol{p} (\ol{\alpha}_o + \mu_k) = \ol{p} (\ol{\alpha}_o + \tau_k - \eta_k) \leq \ol{p} (\ol{\alpha}_o + \tau_k) + \ol{p} (\eta_k) =  \ol{p} (\ol{\alpha}_o + \tau_k).$$
Now, as in the first case, we can extract from $(\mu_k)$ a converging subsequence, yielding thus an element $\ol{\alpha} = \ol{\alpha}_o + \mu$ that satisfies $\ol{d}\ol{\alpha} = \beta$ and $\ol{p}(\ol{\alpha}) = 0$. \\

The fact that $\ol{p}(\ol{\alpha}) = 0$ means that there is a sequence $(\eps_k) \subset B$ such that 
$p(\alpha + \eps_k) \to 0$. Now remembering that $p = p_{r, \{(U_i, \phi_i)\}, \{K_i\}}$, the previous assertion implies that 
\begin{equation}\label{limit} 
\alpha\bigl|_{M_p} = \lim_{k \to \infty} -\eps_k \bigr|_{M_p} \quad \mbox{ \wrt the $C^r$-norm}.
\end{equation}

Now the last step consists in proving that $\alpha$ must therefore be foliated exact on $M_p$. As argued hereafter, this follows mainly from the following two facts. The set of foliated exact forms on $M$ is closed and $M_p$ is compact (as shown in \eref{vanishing-cycle}, the set of foliated forms on the punctured torus endowed with an irrational slope linear foliation is not closed although the same set on the full torus is closed for a non-Liouville slope). \\

Let $K$ be a compact submanifold with boundary of $M$ and let $\Omega^*(K)$ denote the set of restrictions to $K$ of forms on $M$. Equivalently $\Omega^*(K)$ is the quotient of $\Omega^*(M)$ by the closed subspace $\{\alpha ; \alpha|_K = 0\}$. Thus $\Omega^*(K)$ is a nuclear Fr\'echet space as well. It is important to observe that when $K$ is a submanifold with boundary, $\Omega^*(K)$ is defined intrinsically as well and the two definitions coincide. Indeed, taking a partition of unity, the problem reduces to considering the case of a smooth function defined on a relative open subset of the closed upper half space. By this we mean a function that is smooth on the interior of its domain, that admits continuous directional derivatives of all orders for directions that point upwards. It is a well-known result of Whitney\footnote{{\em Whitney's extension theorem} H.~Whitney, Analytic extensions of differentiable functions defined on closed sets, {\em Trans. of the AMS} {\bf 36}, 1934, 63--89.}, also a consequence of the parametric version of the Borel lemma\footnote{M. Golubitsky, V. Guillemin (1974). {\em Stable mappings and their singularities.} Springer-Verlag, Graduate texts in Mathematics: Vol. 14.} that such a smooth function extends across the boundary, that is, extends to a smooth function defined on a genuine open subset of the euclidean space containing the original domain. Now let $T_K\F$ denote the set of vectors at points of $K$ that are tangent to $\F$ and consider the subset $\Omega^*(K, \F_K)$ of forms on $K$ that vanish when evaluated on vector in $T_K\F$. Since the latter is a closed subset, the quotient $\Omega^*(K)/ \Omega^*(K, \F_K)$, also denoted by $\Omega^*(\F_K)$, is a nuclear Fr\'echet space as well. The foliated differential induces a differential
$$d_{\F_K} : \Omega^*(\F_K) \to \Omega^{*+1}(\F_K).$$
Indeed, the de Rham differential induces a differential $d_K$ on $\Omega^*(K)$ that preserves the space $\Omega^*(K,\F_K)$. Then observe that the natural restriction map $r : \Omega^*(\F) \to \Omega^* (\F_K)$ is continuous and surjective. By the open mapping theorem it is thus a homomorphism. Besides, since $\b^*(\F)$ is closed, hence complete, its image through $r$ is complete, hence closed. Finally, the commutativity of the following diagram implies that the image of $d_{\F_K}$ coincides with the set $r (\b^*(\F))$. 

$$\begin{array}{ccc}
\Omega^*(\F) & \stackrel{d_\F}{\longrightarrow} & \Omega^{*+1} (\F) \\
\downarrow r & & \downarrow r \\
\Omega^* (\F_K) & \stackrel{d_{\F_K}}{\longrightarrow} & \Omega^{*+1} (\F_K),
\end{array}$$

Now, the image of $d_{\F_{M_p}}$ being closed, equation (\ref{limit}) implies that the form $\alpha$ is exact on $M_p$. Hence, there exists a foliated form $\gamma$ on $M$ whose restriction to $M_p$ satisfies $\alpha|_{M_p} = d_{\F_{M_p}} \gamma|_{M_p}$. Then the form 
$$\tilde{\alpha} = \alpha - d \gamma.$$
satisfies~:
\begin{enumerate}
\item[-] $\tilde{\alpha}|_{M_p} = 0$, 
\item[-] $d \tilde{\alpha} = \beta$,
\end{enumerate}
that is, we have constructed a form in the kernel of $p$ whose image through $d$ is $\beta$, implying that $\Ker \ul{p} \subset d(\Ker p)$. 
\cqfd

With regards to the properties (DN$_{\rm loc}$) and ($\Omega_{\rm loc}$), their definition can be found in \cite{V}[Definition 2.1 p~815 and Definition 2.4 p~816]. As in the compact case, one proves, using the stability properties established in \cite{V} and the fact that $d_{\F}$ is a SK-homomorphism, that ${\mathcal Z}^*(\F)$ (\rp ${\mathcal B}^*(\F)$) has property ($\Omega_{\rm loc}$) (\rp (DN$_{\rm loc}$)).

\begin{lem}\label{non-compact} If M is non-compact and $H^*(\F)$ is finite-dimensional and Hausdorff then $d_{\F}$ admits a continuous linear right inverse.
\end{lem}

\Pf Let us introduce some notation~: $E_k = \Omega^k(\F)$ and $T_k = d_\F|_{\Omega^k(\F)}$, with $E_{-1} = \{0\}$ and $T_{-1} = 0$. We are thus considering a non-exact sequence
$$E_{k-1} \stackrel{T_{k-1}}{\longrightarrow} E_k  \stackrel{T_k}{\longrightarrow} E_{k+1}  \stackrel{T_{k+1}}{\longrightarrow} E_{k+2}$$
such that 
\begin{enumerate}
\item[-] $T_k \circ T_{k-1} = 0$ and $T_{k+1} \circ T_k = 0$,
\item[-] $\Ker T_{k} / \im T_{k-1} \stackrel{{\rm not}}{=} F_k$ and $\Ker T_{k+1} / \im T_k \stackrel{{\rm not}}{=} F_{k+1}$ are finite-dimensional.
\end{enumerate}  
Because the map $T_{k+1}$ is a SK-homomorphism, the space $\Ker T_{k+1}$ is a SK-subspace (\cite{V}[Definition 1.4, p~814]). Therefore $\Ker T_{k+1}$ inherits the property (DN$_{\rm loc}$) from $E_{k+1}$ (cf.~\cite{V}[Lemma 2.3 p~816]). Similarly, the space $\im T_{k-1}$ has the property ($\Omega_{\rm loc}$) because $T_{k-1}$ is a SK-homomorphism (cf.~\cite{V}[Lemma 2.6 p~817]).\\

Now, in order to prove that $\im T_k = \b^{k+1}(\F)$ has property (DN$_{\rm loc}$), it is sufficient, again according to \cite{V}[Lemma 2.3 p~816] to prove that $\im T_k$ is a SK-subspace of $\Ker T_{k+1}$, that is, to prove that the projection $\Ker T_{k+1} \to \Ker T_{k+1} / \im T_k = F_{k+1}$ is a 
SK-homomorphism. This follows directly from the fact that $F_{k+1}$, being finite-dimensional and Hausdorff, carries a compatible norm (or from the fact that the projection $\Ker T_{k+1} \to F_{k+1}$ has a continuous right inverse). \\

To prove that $\Ker T_k = \ZE^k(\F)$ has property ($\Omega_{\rm loc}$), we use the fact that $\Ker T_k \cong \im T_{k-1} \oplus F_k$ is the image of the map 
$$T_{k-1} \oplus \id : E_{k-1} \oplus F_k \to E_k \oplus F_k,$$ 
which, as we verify hereafter, is a SK-homomorphism. Let $p$ be a continuous seminorm on $E_{k-1} \oplus F_k$. It induces a continuous seminorm $p_o$ on $E_{k-1}$ by restriction and hence a continuous seminorm $q$ on $E_k$ such that 
$$T_{k-1} (\Ker p_o) \supset \Ker q \cap \im T_{k-1}.$$
Let $\ol{q}$ be an extension of $q$ to $E_k \oplus F_k$ of the type $q + ||\cdot||$, where $||\cdot||$ is a norm on $F_k$. Then $\Ker \ol{q} = \Ker q$ and hence
$$(T_{k-1} \oplus \id) (\Ker p) \supset (T_{k-1} \oplus \id) (\Ker p_o \oplus \{0\}) \supset \Ker\ol{q} \cap \im (T_{k-1} \oplus \id).$$
Moreover, the space $E_{k-1} \oplus F_k$ is ($\Omega_{\rm loc}$) for $k \geq 1$ because $s^\N \oplus F_k \cong s^\N$ and for $k=0$ because a finite-dimensional Hausdorff space trivially satisfies the property ($\Omega_{\rm loc}$).
\cqfd

\noindent
{\em End of the proof of \pref{finiteK}.}

\begin{rmk}
In the case of a one-leaf foliation the fact that the de Rham differential admits a right inverse is easily seen in the compact case either using Hodge theory when the manifold is orientable or using \v{C}ech cohomology of a finite good cover in general. In the non-compact case, the approach used here requires to prove that $d$ is a SK-homomorphism, that the space ${\mathcal B}^*(M)$ of exact forms has property (DN$_{\rm loc}$) and that the space ${\mathcal Z}^*(M)$ of closed forms has property ($\Omega_{\rm loc}$). I do not know how to prove that last two assumptions when the rank of the cohomology of $M$ is infinite. 
That $d$ is a SK-homomorphism can be approached as follows~: it is sufficient to have an exhaustion of the manifold $M$ by compact sets $K_i \subset K_{i+1}$ such that the various restriction maps $H^*(M) \to H^*(K_i)$ are surjective. In other words, the sets $K_i$ do not carry ``useless" cohomology classes, that is, cohomology classes that do not come from the ambient manifold. Such and exhaustion is called hereafter an {\em economical exhaustion}. A proper Morse function for which the Morse number equals the type number (that is, the number of critical points coincides with the rank of the cohomology) yields such a decomposition. Such a Morse function exists on any simply-connected manifolds of dimension $>5$ with torsionless homology\footnote{cf.~John Franks, The periodic structure of nonsingular Morse-Smale flows. {\em Comment. Math. Helv.} {\bf 53}, (1978), no. 2, 279--294.}. However, it seems much easier to construct directly an economical exhaustion and I believe it exists on any open manifold.
\end{rmk}

\section{Counterexample}\label{ex}

One would like to understand what happens when neither of the situations encountered above occurs. Let $(M, \F)$ and $(N,\G)$ be two foliated manifolds. Suppose that both foliated cohomologies are infinite-dimensional with one of them non-Hausdorff. Then the tensor product $H^*(\F) \otimes H^*(\G)$ cannot be completed. There is nevertheless a case where an alternative to completion could be proposed, that is, when one of the foliations, say $\G$, is a foliation by points $\G = \F_N$. Then $H^*(\G)$ coincides with $C^\infty (N)$ and one is tempted to replace $H^*(\F) \hot C^\infty (N)$, which does not make sense here, by $C^\infty(N, H^*(\F))$, since these two spaces coincide when $H^*(\F)$ is Hausdorff (cf.~\cite{T}[Theorem 44.1 p~449]). It is therefore natural to wonder whether the trivial map 
$H^*(\F \times \F_N) \to C^\infty(N, H^*(\F_N))$ yields an isomorphism or not~:
\begin{equation}
H^*(\F \times \F_N) \stackrel{?}{\cong} C^\infty(N, H^*(\F)).
\end{equation}

The answer is negative. Indeed, the torus ${\mathbb T}^2$ endowed with a Liouville foliation (see \eref{torus}) supports a smooth family of foliated exact forms which is not the coboundary of any smooth, nor even continuous, family of forms. This smooth family represents thus both the zero element in $C^\infty(N, H^*(\F))$ and a non-zero element in  $H^*(\F \times \F_N)$. I do not see any obvious theoretical reason for such a family to exist; both the space $H^*(\F \times \F_N)$ and the space $C^\infty(N, H^*(\F))$ are non-Hausdorff; somehow $H^*(\F \times \F_N)$ is ``more separated" than $C^\infty(N, H^*(\F))$. \\ 

Let $x, y$ denote standard coordinates on the torus $\torus^2$. The leaves of the foliation $\F_\alpha$ are the orbits of the vector field $X = \partial_x + \alpha \partial_y$. Any foliated $1$-form is automatically closed and can be written $f \ol{dx}$, with $f$ in $C^\infty(\torus^2)$ and $\ol{dx}$ the image of the closed form $dx \in \Omega^1({\mathbb T}^2)$ in $\Omega^1(\F_\alpha)$. It is exact when $f \ol{dx} = \ol{dg}$, for some $g$ in $C^\infty(\torus^2)$, which is equivalent to $f = Xg$. Besides, we may consider the Fourier expansions of the functions $f$ and~$g$~:

$$f = \sum_{m,n \in \Z} f_{m,n} e^{2 \pi i (m x + n y)}  \quad \mbox{and} \quad g = \sum_{m,n \in \Z} g_{m,n} e^{2 \pi i (m x + n y)}.$$
The equation $f = Xg$ is equivalent to the sequence of equations 

$$f_{m,n} = 2 \pi i (m + \alpha n) \; g_{m,n}, \quad m, n \in \Z,$$
which of course implies $f_{0,0} = 0$. 

\begin{lem}\label{smooth}
Suppose that $(f_t)_{t \in \R}$ is a family of functions on $\torus^2$. It is a smooth family of smooth functions if and only if each function $t \mapsto (f_t)_{m,n}$ is smooth and for all compact interval $I$ in $\R$, and integers $a \geq 0$, $j \geq 1$, there exists a constant $c = c(I,a,j)$ such that 

$$\sup_{t \in I} \big| \partial^a_t (f_t)_{m,n} \big| \leq \frac{c}{(|m|+|n|)^j}.$$
\end{lem}
\noindent
To see the necessity of this condition it suffices to combine part integration in order to get rid of the derivatives \wrt $x$ and $y$ with the fact that the Fourier coefficient $(\partial_x^k \partial_y^l \partial_t^a f_t)_{m,n}$ is bounded by a constant depending only on $I, k, l$ and $a$. \\

With these preliminaries in mind, we are ready to construct a family $f_t$ of functions on $\torus^2$ with the following properties~:

\begin{enumerate}
\item[(i)] $f_t$ is a smooth family of smooth functions,
\item[(ii)] for each value of the parameter $t$, there is a smooth solution to $f_t = Xg_t$,
\item[(iii)] there no smooth --- nor even continuous --- family of smooth functions $g_t$ solving $f_t = 
Xg_t$.
\end{enumerate}
Since $\alpha$ is a Liouville number, for each integer $ p > 1$, there exists a pair of integers $(m_p, n_p)$ such that 
$$|m_ p + \alpha n_ p| \leq \frac{1}{(|m_ p|+|n_ p|)^ p}.$$
\Wlog assume that $(m_p, n_p) \neq (m_q, n_q)$ for $p \neq q$ and that $n_p \geq p$. Now define

$$(f_t)_{m,n} = \left\{ 
\begin{array}{ll}
(m_ p + \alpha n_ p) (|m_p|+|n_p|)\rho(s_ p(t - \frac{1}{ p})) & \mbox{if}  \quad (m,n) = (m_ p,n_ p), \\
0 & \mbox{otherwise} 
\end{array}
\right.$$
where $\rho$ is a bump function supported in the interval $[-1, 1]$ that achieves its maximum value $1$ at $0$ and where $s_p = p (p + 1)$. The function $\rho(s_p(t - \frac{1}{p}))$ has its support contained in $[\frac{1}{p} - \frac{1}{2p(p+1)}, \frac{1}{p} + \frac{1}{2p(p+1)}]$. \\

Let us verify that the $(f_t)_{m,n}$'s are the Fourier coefficients of a family $f_t$ enjoying the properties (i), (ii) and (iii). \\

\noindent
(i) For smoothness of $f_t$ we use the criterion described in \lref{smooth}.
$$\begin{array}{ccl}
|\partial_t^a (f_t)_{m_p,n_p}| & \leq & \displaystyle{\big| (m_p + \alpha n_p) \big| (|m_p|+|n_p|) \sup_{t \in I} |\partial_t^a \rho(t)| |s_p|^a} \\
& \leq & \displaystyle{\frac{c_a |s_p|^a}{(|m_p|+|n_p|)^{p - 1}}} \\
& \leq & \displaystyle{\frac{c'_a}{(|m_p|+|n_p|)^{p - 1 - 2a}}} \\
& \leq & \displaystyle{\frac{c''_{a,j}}{(|m_p|+|n_p|)^j}}
\end{array}$$

The before-last inequality follows from the fact that $s_p$ is a polynomial of degree $2$ in $p$  and the assumption $n_p \geq p$, while the last inequality is a consequence of the fact that $p - 1 - 2a \to \infty$ when $p \to \infty$. \\

\noindent
(ii) The coefficients $$(g_t)_{m,n} = (f_t)_{m,n} /(m + \alpha n)$$ define a smooth function for each value of $t$. Indeed, for a fixed $t_0$,

$$(g_{t_0})_{m_p,n_p} = (|m_p|+|n_p|) \rho(s_ p(t_0 - \frac{1}{ p})) = 0$$
for all $p$'s except perhaps one since the supports of the various functions $\rho(s_ p(t_0 - \frac{1}{ p}))$ are disjoints. The Fourier series of the function $g_{t_0}$ has thus only one term. \\

\noindent
(iii) The function $g_t$ is not smooth, nor even continuous, near $t=0$. Indeed, the coefficients $(g_t)_{m,n}$ are not uniformly bounded on any interval $I$ around $0$~: 

$$\sup_{t \in I} |(g_t)_{m_p,n_p}| = (|m_p|+|n_p|) \; \sup_{t \in I} |\rho(s_ p(t_0 - \frac{1}{ p}))| = (|m_p|+|n_p|)$$
as soon as $p$ is sufficiently large for $\frac{1}{p}$ to belong to $I$.

\end{document}